# Grothendieck's Double Limit Theorem and Model Theory

Karim Khanaki

Arak University of Technology, &

Institute for Research in Fundamental Sciences (IPM)

Abstract: This is an expository paper in Persian on Grothendieck's double limit theorem and its connection with (neo-)stability project. We review recent results/observations and also discuss historical and philosophical issues.

## ۱. گروتندیک، قضیهٔ حدِ دوگانه، ۱۹۵۲

در کنگرهٔ بین‌المللی ریاضی‌دانان در مسکو (۱۹۶۶)، الکساندر گروتندیک[۱]، اساساً برای دستاوردهایش در هندسه جبری، برندهٔ مدال فیلدز شد؛ در حالی‌که برخی از بزرگان آنالیز از جمله دیودونه[۲]، معتقد بودند که دستاوردهای او در فضاهای برداری توپولوژیک (۱۹۵۰–۱۹۵۵) برای کسب این افتخار کافی‌ست. کار او در این دوره و تاثیر آن در سال‌های بعد را با کارهای باناخ[۳] مقایسه کرده‌اند.[۴] هدف این مقاله مطالعهٔ یکی از دستاوردهای مهم گروتندیک در این دورهٔ زمانی است که به‌طور غیرمنتظره‌ای در مرکزِ نظریهٔ مدل‌های مدرن قرار می‌گیرد. در این فصل، این دستاورد گروتندیک، موسوم به قضیهٔ حدِ دوگانه[۵]، ارائه می‌گردد و در فصل‌های بعد ارتباط آن با نظریهٔ مدل‌ها مورد مطالعه قرار می‌گیرد.

در ادامهٔ این مقاله، $C(X)$ فضای باناخِ همهٔ توابع پیوستهٔ حقیقی-مقدار روی یک فضای توپولوژیک فشردهٔ $X$، با نُرم یکنواخت[۶]، است. برای یک زیرمجموعه $A$ از $\mathbb{R}^X$، توپولوژی همگرایی نقطه‌ای[۷] روی $A$ عبارت است از توپولوژی زیرفضایی به‌ارث‌رسیده از فضای حاصل‌ضربی $\prod_{x\in X}\mathbb{R}$. یادآوری این نکته مفید است که برای یک زیرمجموعهٔ دلخواه $A$ از $C(X)$، بستار $A$ در توپولوژی همگرایی نقطه‌ای می‌تواند شامل توابعی باشد که حتی اندازه‌پذیر بورل[۸] نباشند، چه رسد به آنکه پیوسته باشند.

---

[1] Alexander Grothendieck (1928-2014)

[2] Jean Dieudonné (1906-1992)

[3] Stefan Banach (1892-1945)

[۴] مراجعه شود به [۱۳]، صفحه ۶۴۴.

[5] Grothendieck's double limit theorem (DLT)

[6] Uniform norm/sup-norm

[7] Topology of pointwise convergence

[8] Borel measurable

قضیه ۱.۱ (قضیهٔ حدِ دوگانهٔ گروتندیک). فرض کنید $X$ یک فضای توپولوژیک فشرده و $X_0$ زیرمجموعه چگالی از آن باشد. آنگاه برای هر زیرمجموعه کران‌دار $A$ از $C(X)$ دو گزارهٔ زیرین هم‌ارز هستند:

(۱) بستار $A$ در توپولوژی همگرایی نقطه‌ای همچنان زیرمجموعه‌ای از $C(X)$ است.

(۲) هرگاه $(f_n)$ و $(x_m)$ به‌ترتیب دنباله‌هایی در $A$ و $X_0$ باشند، آنگاه

$$\lim_n \lim_m f_n(x_m) = \lim_m \lim_n f_n(x_m),$$

با این شرط که هر دو حد موجود باشند.

اجازه دهید پیش از ارائه اثبات این قضیه چند نکته را بیان کنیم:

ملاحظه ۲.۱.

(۱) قضیهٔ فوق حالت خاصی از قضیه ۶ در [۷] می‌باشد. در واقع، اگر $X$ یک فضای توپولوژیک دلخواه باشد و بستار $A$ در توپولوژی ضعیف[1] درنظر گرفته شود، هم‌ارزی فوق همچنان صادق است. توجه داریم که اگر $X$ فشرده باشد، توپولوژی همگرایی نقطه‌ای و توپولوژی همگرایی ضعیف، روی زیرمجموعه‌های کران‌دار از $C(X)$، یکسان هستند (مراجعه شود به [۵]، بخش ۴۶۲).

(۲) برای اهداف این مقاله، این حالتِ خاصِ قضیهٔ گروتندیک (یعنی قضیه ۱.۱) نیز کافی خواهد بود. در واقع، در فصل ۳ خواهیم دید که $X$ به‌صورت فضای تایپ‌ها[2] و $A$ خانواده‌ای از فرمول‌ها تعبیر خواهند شد.

(۳) گزاره (۲) در قضیهٔ فوق «ویژگی حدِ دوگانه[3]» نامیده می‌شود.

اثبات قضیه ۱.۱. اثبات (۱) $\Longleftarrow$ (۲) سخت نیست، و به خواننده واگذار می‌شود.

(۲) $\Longleftarrow$ (۱): برای به‌دست آوردنِ یک تناقض، فرض کنید که $f$ یک تابع در بستار $A$ باشد که در نقطه $x \in X$ پیوسته نیست. بنابراین، یک همسایگی $U$ از $f(x)$ وجود دارد به‌طوری‌که هر همسایگی از $x$ شامل نقطه‌ای مانند $y$ است به‌طوری‌که $f(y) \notin U$. تابع $f_1$ از $A$ را درنظر بگیرید؛ آنگاه یک نقطه $x_1$ وجود دارد به‌طوری‌که $|f_1(x_1)-f_1(x)| < 1$ و $f(x_1) \notin U$. تابع $f_2$ در $A$ را طوری درنظر بگیرید که $|f_2(x_1)-f(x_1)| < 1$ و $|f_2(x)-f(x)| < 1$. اکنون نقطهٔ $x_2$ را طوری انتخاب کنید که $|f_i(x_2)-f_i(x)| < \frac{1}{2}$ $(i=1,2)$ و $f(x_2) \notin U$. تابع $f_3$ را طوری انتخاب کنید که $|f_3(x_i)-f(x_i)| < \frac{1}{2}$ و $|f_3(x)-f(x)| < \frac{1}{2}$. با ادامه این فرآیند، دنباله‌های $(x_m)$ و $(f_n)$ را می‌یابیم به‌طوری‌که برای هر $n$، $|f_i(x_n) - f_i(x)| < \frac{1}{n}$

[1] weak topology
[2] space of types
[3] Double Limit Property (DLP)

$|f_{n+1}(x_j) - f(x_j)| < \frac{1}{n}$ ،$(j = 1, \ldots, n)$ $|f_{n+1}(x) - f(x)| < \frac{1}{n}$ ،$(i = 1, \ldots, n)$

و $f(x_n) \notin U$. اکنون، به‌آسانی می‌توان بررسی نمود که گزاره (۲) در قضیه ۱.۱ نقض می‌شود.[۱] □

نکات ذیل برای خواننده‌ای که اثبات فوق را ساده‌می‌انگارد حائز اهمیت است: یکم، همانطور که قبلاً بیان شد، قضیهٔ فوق حالتِ خاصِ قضیهٔ گروتندیک است، و اهمیت و کاربرد آن در آنالیز و نظریه مدل‌ها در ادامه روشن خواهد شد. دوم و مهمتر از منظر روش‌شناختی اینکه، برای یک نظریه‌مدل‌دان که با مفاهیمِ قدرتمندِ توپولوژی نمی‌اندیشد، ارائه چنین اثباتی در زبان نظریهٔ مدل‌ها نمی تواند به‌سادگی و وضوح اثبات فوق باشد.

پیش از بررسی ارتباط قضیهٔ گروتندیک با نظریهٔ مدل‌ها (در فصل‌های بعد)، شایسته است برای مشاهدهٔ کاربردهایی از این قضیه در آنالیز تابعی خواننده را به [۲۴] ارجاع دهیم. برای مشاهدهٔ کاربرد بسیار زیبایی از خاصیت حد دوگانه در نظریه فضاهای باناخ، خواننده را به ملاحظه ۲.۴(۴) در ذیل ارجاع می‌دهیم.

مباحث ارائه شده در این مقاله به شرح ذیل می‌باشد. در فصل بعد قضیهٔ بنیادین نظریهٔ پایداری[۲] (قضیه ۲.۲) ارائه می‌گردد. در فصل ۳ نشان می‌دهیم که قضیهٔ بنیادین نظریهٔ پایداری نتیجه‌ای از قضیهٔ گروتندیک است. در فصل ۴ معنای نظریه مدلی قضیهٔ گروتندیک را روشن می‌کنیم. در فصل ۵ بهبودی از قضیهٔ گروتندیک برای تایپ‌های به‌طور ژِنِریک پایدار[۳] ارائه می‌دهیم. در فصل آخر در مورد پرسش‌های باقیمانده و کارهای آتی بحث خواهیم کرد.

## ۲. شِلاح، قضیهٔ بنیادینِ نظریهٔ پایداری، ۱۹۷۱

حدوداً دو دهه بعد از قضیهٔ گروتندیک، در [۲۱]، شِلاح[۴] قضیهٔ بنیادین نظریهٔ پایداری (قضیه ۲.۲) را ثابت کرد. خواهیم دید که قضیهٔ اخیر نتیجه‌ای از قضیه ۱.۱ است.

فرض می‌کنیم که خواننده با مفاهیم مقدماتی نظریهٔ مدل‌ها که در کتاب‌های منطق ریاضی موجوداند آشنایی دارد، از جمله: زبان، فرمول، تئوری، مدل، محقق شدن[۵] (یا ارضاء شدن). البته، پیش از ارائهٔ قضیهٔ بنیادین نظریهٔ پایداری می‌بایست چند مفهوم دیگر را معرفی کنیم.

در ادامهٔ این مقاله، $L$ یک زبان منطق مرتبهٔ اول (کلاسیک)، $\varphi(x,y)$ یک $L$-فرمول، $T$ یک $L$-تئوری و $M$ یک مدل از $T$ است.

**تعریف ۱.۲.**

---

[۱] توجه داریم که چون $X_\circ$ چگال است، می‌توانیم دنبالهٔ $(x_n)$ را در $X_\circ$ انتخاب کنیم.

[2] Fundamental Theorem of Stability (FTS)

[3] generically stable types

[4] Saharon Shelah (1945-)

[5] satisfaction

(۱) گوییم $\varphi(x,y)$ در $M$ ویژگی ترتیب[1] دارد (و به‌طور مُخفّف می‌نویسیم $OP$ دارد)، هرگاه دنباله‌های $(a_n)$ و $(b_n)$ در $M$ موجود باشند به‌طوری‌که $\varphi(a_m,a_n)$ در $M$ درست است اگر و تنها اگر $m<n$. اگر $\varphi(x,y)$ در $M$ ویژگی ترتیب نداشته باشد، گوییم $\varphi(x,y)$ در $M$ پایدار است.

(۲) گوییم $\varphi(x,y)$ در تئوری $T$ پایدار است هرگاه برای هر مدل $M$ از $T$، $\varphi(x,y)$ در $M$ پایدار باشد. به چنین فرمولی اصطلاحاً «فرمولِ پایدار[2]» گویند.

در ادامه تعریفِ $\varphi$-تایپ ارائه می‌شود که مشابه تعریف تایپ‌ها در نظریه‌مدل‌ها است با این تفاوت که در تعریف آن فقط از فرمول‌هایی به‌صورت $\varphi(x,b)$ یا $\neg\varphi(x,b)$ استفاده می‌شود. در واقع، این مفهوم به نوعی موضعی‌شدهٔ مفهوم تایپ است، به این معنا که تعریف آن فقط وابسته به فرمول $\varphi$ است و نه همهٔ فرمول‌های زبان.

یک $\varphi$-تایپ روی $M$ مجموعه‌ای سازگار از فرمول‌های $\varphi(x,b)$ یا $\neg\varphi(x,b)$ است که $b\in M$.[3] برای $a\in M$، $\varphi$-تایپِ $a$ روی $M$، که با $tp_\varphi(a/M)$ نشان داده می‌شود، مجموعهٔ همهٔ فرمول‌های $\varphi(x,b)$ یا $\neg\varphi(x,b)$ ($b\in M$) است که در $M$ توسط $a$ محقق می‌شوند. یک $\varphi$-تایپ $p$ را کامل گوییم اگر هیچ $\varphi$-تایپ دیگری شامل $p$ وجود نداشته باشد. (به‌وضوح، هر تایپ $tp_\varphi(a/M)$ کامل است.) مجموعهٔ همهٔ $\varphi$-تایپ‌های کامل روی $M$ را با $S_\varphi(M)$ نشان می‌دهیم.

یک تایپ $p$ در $S_\varphi(M)$ را تعریف‌پذیر روی $M$[4] گوییم هرگاه یک فرمول $\psi(y)$ با پارامترهایی در $M$ وجود داشته باشد به‌طوری‌که برای هر $b\in M$، داریم $\varphi(x,b)\in p$ اگر و فقط اگر $\psi(b)$ درست باشد. (تعریف‌پذیری مفهومی بسیار مهم در نظریهٔ مدل‌ها است که به‌طور مستقل در دههٔ هفتاد توسط شلاح و گایفمن[5] معرفی شده است. در فصل ۳ خواهیم دید که این مفهوم متناظر با مفهوم «پیوستگی» در توپولوژی است، لذا اهمیت آن کاملاً طبیعی است.)

اکنون آماده‌ایم که قضیهٔ بنیادین را بیان کنیم.

قضیه ۲.۲ (قضیهٔ بنیادینِ نظریهٔ پایداری شِلاح). فرض کنید $\varphi(x,y)$ یک فرمول و $T$ یک تئوری باشند. آنگاه سه گزارهٔ زیرین هم‌ارز هستند:

(۱) $\varphi(x,y)$ در تئوری $T$ پایدار است.

---

[1] The Order Property (OP)

[2] stable formula

[3] به زبانی ساده، یک $\varphi$-تایپ را می‌توان مجموعه‌ای از معادلات (یا دستگاه معادلات) تصور کرد که جوابی مشترک برای همهٔ آنها وجود دارد، البته در یک مدل بزرگ.

[4] definable over $M$

[5] Haim Gaifman (1934-)

(۲) برای هر مدل $M$ از $T$، هر $\varphi$-تایپ روی $M$ تعریف‌پذیر است.[۱]

(۳) برای هر مدل $M$ از $T$، تعداد $\varphi$-تایپ‌های روی $M$ حداکثر $|M|$ است.[۲]

اثبات (۳) ⇐ (۱) استاندارد است و در بیشتر کتاب‌های نظریهٔ مدل‌ها یافت می‌شود. (۲) ⇐ (۳) از این حقیقت نتیجه می‌شود که تعداد فرمول‌هایی که به‌صورت ترکیبات بولی متناهی از فرمول‌های $\varphi(a,y)$ (که $a \in M$) هستند، حداکثر $|M|$ است. در فصل بعد نشان خواهیم داد که چگونه قسمت اصلی قضیهٔ بنیادین، یعنی (۱) ⇐ (۲)، نتیجه‌ای از قضیه ۱.۱ است. برای این کار نیاز به ترجمهٔ قضیهٔ گروتندیک به زبان نظریهٔ مدل‌ها داریم.

ملاحظه ۳.۲.

(۱) اثبات شلاح برای (۱) ⇐ (۲) کاملاً با اثباتی که در فصل بعد (با استفاده از قضیهٔ گروتندیک) ارائه می‌گردد متفاوت است. درواقع، اثبات‌های کلاسیک محتوای «ترکیبیاتی[۳]» دارند.

(۲) پس از قضیهٔ بسیار مهم مورلی[۴] در ۱۹۶۵ موسوم به قضیه جازمیت[۵]، مفهوم تایپ و تئوری‌های پایدار در کانون نظریهٔ مدل‌های نوین قرار گرفت. تا کنون، مثال‌های مهمی از تئوری‌های پایدار مورد مطالعه قرار گرفته‌اند و کاربردهای بسیاری از نظریهٔ پایداری در حوزه‌های دیگر ریاضیات یافت شده‌اند. به عنوان نمونه، هُروشوفسکی[۶] در دههٔ نود اثباتی از حدس هندسی موردِل-لَن[۷] در همهٔ مشخصه‌ها ارائه داد.

(۳) بسط نتایج به‌دست آمده در تئوری‌های پایدار[۸] به تئوری‌های ناپایدار از اهداف مهم نظریه مدل‌دان‌ها در قرن بیست و یکم است. یکی از این نتایج در فصل ۵ ارائه خواهد شد.

## ۳. بِن‌یاکوف، ترجمهٔ قضیهٔ گروتندیک به زبان نظریهٔ مدل‌ها، ۲۰۱۵

بن‌یاکوف[۹] [۲] نقل می‌کند:

> هنگامی که از حضّار در جلسه‌ای در کلکته (۲۰۱۳) پرسش کردم «چه کسی، و در چه زمانی، برای اولین‌بار مفهوم فرمول پایدار را تعریف کرد»، در مخالفت با جواب آنها پاسخی دادم که شبیه یک شوخی بود: «گروتندیک، در دههٔ پنجاه.»[۱۰]

---

[۱] علاوه بر این، فرمول‌هایی که $\varphi$-تایپ‌ها را تعریف می‌کنند به‌صورت ترکیبات بولی متناهی از فرمول‌های $\varphi(a,y)$ هستند که $a \in M$.

[۲] گزاره (۳) دلیل نام‌گذاری تئوری‌های پایدار می‌باشد، زیرا در چنین تئوری‌هایی تعداد تایپ‌ها پایدار می‌ماند.

[3] combinatorial

[4] Michael Darwin Morley (1930-2020)

[5] Morley's Categoricity Theorem

[6] Ehud Hrushovski (1959-)

[7] Mordell–Lang conjecture

[8] unstable theories

[9] Itaï Ben Yaacov

[۱۰] برای نظریه مدل‌دان‌ها پاسخ معقول و شناخته شده می‌بایست «شلاح، در دههٔ هفتاد» باشد.

بیان دقیق‌تر این است که قضیهٔ بنیادین نظریهٔ پایداری، یعنی معادل بودنِ «تعریف‌پذیری تایپ‌ها» و «نداشتن خاصیت ترتیب»، محتوای قضیهٔ گروتندیک است، هرچند به عنوان یک دستاورد در نظریهٔ مدل‌ها، این قضیه توسط شلاح در [۲۱] ثابت شده است.

در ادامه نشان خواهیم داد که بن‌یاکوف چگونه محتوای قضیهٔ گروتندیک را به زبان نظریهٔ مدل‌ها ترجمه می‌کند و با استفاده از آن اثباتی از قضیهٔ بنیادین نظریهٔ پایداری ارائه می‌دهد.

تعریف می‌کنیم $\tilde{\varphi}(y,x) := \varphi(x,y)$ (یعنی نقش $x$ و $y$ عوض شده است). مشابه فصل قبل، مجموعهٔ همهٔ $\tilde{\varphi}$-تایپ‌های کامل روی $M$ را با $S_{\tilde{\varphi}}(M)$ نشان می‌دهیم.

برای هر $a \in M$ تابع $f_a : S_{\tilde{\varphi}}(M) \to \{۰,۱\}$ را به این‌صورت تعریف می‌کنیم که $f_a(q) = ۱$ اگر برای یک تحقق $b$ از $q$، $\varphi(a,b)$ درست باشد، و $f_a(q) = ۰$ در غیر این‌صورت. به‌وضوح، اطلاعات تایپ $tp_\varphi(a/M)$ توسط این تابع کُد می‌شود. می‌توان $S_{\tilde{\varphi}}(M)$ را به توپولوژی مجهز نمود که با آن همهٔ توابع $f_a$ ($a \in M$) پیوسته باشند. دو نکته اساسی این است که: (۱) هر تایپ در $S_\varphi(M)$ توسط تابعی در بستار مجموعه $A = \{f_a : a \in M\}$ (در توپولوژی همگرایی نقطه‌ای) کُد می‌شود، و (۲) یک تایپ در $S_\varphi(M)$ تعریف‌پذیر است اگر و تنها اگر تابع متناظر آن پیوسته باشد. (بررسی این دو نکته به خواننده واگذار می‌شود. در مورد (۱)، توجه کنید که تایپ‌هایی که به‌صورت $tp_\varphi(a/M)$ ($a \in M$) هستند در $S_\varphi(M)$ چگال‌اند. برای (۲)، از آنجا که توپولوژیِ $S_{\tilde{\varphi}}(M)$ تماماً ناهمبند است، یک تابع $\{۰,۱\}$-مقداری روی آن پیوسته است اگر و تنها اگر توسط یک «فرمول» با پارامترهایی در $M$ تعریف شود.)

همچنین، اگر تعریف کنیم $\varphi(a,b) := ۱$ هرگاه $\varphi(a,b)$ درست باشد، و $\varphi(a,b) := ۰$ در غیر این‌صورت، آنگاه به‌آسانی می‌توان بررسی نمود که پایدار بودن $\varphi(x,y)$ روی $M$، هم‌ارز است با اینکه برای هر دو دنبالهٔ $(a_m)$ و $(b_m)$ در $M$ داشته باشیم

$$\lim_n \lim_m \varphi(a_m, b_n) = \lim_m \lim_n \varphi(a_m, b_n),$$

با این شرط که هر دو حد موجود باشند.

اکنون آماده‌ایم که اثبات قضیهٔ بنیادینِ نظریهٔ پایداری (۲.۲) را کامل کنیم.

اثبات (۱)⇐(۲) از قضیه ۲.۲. فرض کنید $X$ فضای $S_{\tilde{\varphi}}(M)$ و $X_۰$ مجموعهٔ همهٔ تایپ‌های $tp(a/M)$ (برای $a \in M$) باشد. توجه داریم که $X$ فشرده است و $X_۰$ در $X$ چگال است. همچنین، مجموعه $A = \{f_a : a \in M\}$ زیرمجموعه‌ای کران‌دار از $C(X)$ است. بنابراین، با توجه به پایدار بودن $\varphi(x,y)$ روی $M$، قضیه ۱.۱ نتیجه می‌دهد که هر تابع در بستار $A$ پیوسته است. از آنجا که هر تایپ در $S_\varphi(M)$ توسط تابعی در بستار

مجموعه $A$ کُد می‌شود، و یک تایپ در $S_\varphi(M)$ تعریف‌پذیر است اگر و تنها اگر تابع متناظر آن پیوسته باشد، اثبات کامل است.[۱] $\square$

در ادامه کنکاشی عمیق‌تر در معنای نظریه مدلی قضیهٔ گروتندیک خواهیم داشت.

## ۴. پیلِی، معنای نظریه مدلی قضیهٔ گروتندیک، ۲۰۱۶

در فاصلهٔ کوتاهی از مشاهدات بن‌یاکوف، پیلی[۲] در مقالهٔ [۱۶] تصویر روشن‌تری از معنای نظریه مدلی قضیهٔ گروتندیک ارائه داد. درواقع، این قضیه هم‌ارز با تعریف‌پذیری تایپ‌های جهانی[۳] است، که در ادامه شرح داده خواهد شد.

یک مدل که همهٔ تایپ‌ها (روی مدل‌هایی به اندازهٔ کافی بزرگ) در آن محقق می‌شوند را مدل بزرگ[۴] گوییم، و آنرا با $\mathbb{M}$ نشان می‌دهیم.

در تعریف $\varphi$-تایپ کامل در بخش قبل، اگر مدل $M$ را با مدل بزرگ $\mathbb{M}$ جایگزین کنیم، آنگاه یک $\varphi$-تایپ جهانی خواهیم داشت. یک $\varphi$-تایپ جهانی $p(x)$ را تعریف‌پذیر روی $M$ گوییم هرگاه یک فرمول $\psi(y)$ (با پارامترهایی در $M$) وجود داشته باشد به‌طوری‌که برای هر $b \in \mathbb{M}$، داریم $p \ni \varphi(x,b)$ اگر و فقط اگر $\psi(b)$ درست باشد. یک $\varphi$-تایپ جهانی $p$ را متناهیاً محقق شده در $M$[۵] گوییم، هرگاه هر فرمول $\phi(x)$ که متعلق به $p$ باشد در $M$ محقق شود.

با نمادگذاری بخش قبل، فرض کنید $X = S_{\tilde{\varphi}}(M)$ و $A = \{f_a : a \in M\}$. در ملاحظه ۱.۲ در [۱۶]، پیلی نشان داد که تناظری یک‌به‌یک بین همهٔ توابع در بستار مجموعه $A$ (در توپولوژی همگرایی نقطه‌ای) و مجموعهٔ همهٔ $\varphi$-تایپ‌های جهانی که متناهیاً در $M$ محقق می‌شوند وجود دارد.

مشاهده اخیر، این امکان را به پیلی می‌داد که تصویر کاملی از معنای نظریه مدلی قضیهٔ گروتندیک ارائه دهد:

**قضیه ۱.۴** (معنای نظریه مدلی قضیهٔ گروتندیک). فرض کنید $\varphi(x,y)$ یک فرمول و $M$ مدلی از تئوری $T$ باشد. آنگاه دو گزارهٔ زیرین هم‌ارز هستند:

(۱) $\varphi(x,y)$ در مدل $M$ پایدار است.

---

[۱] همچنین، به‌آسانی می‌توان نشان داد که تابع پیوسته‌ای که هر تایپ را تعریف می‌کند متناظر با یک فرمول به‌صورت ترکیب بولی متناهی از فرمول‌های $\varphi(a,y)$ است که $a \in M$.

[2] Anand Pillay (1951-)

[3] global types

[۴] تعریف دقیق مدل‌های بزرگ کاملاً فنّی است و خارج از بحث این مقاله می‌باشد. خواننده می‌تواند به مبحث monster model در کتاب‌های نظریهٔ مدل‌ها مراجعه کند.

[5] finitely satisfiable in $M$

(۲) هر تایپ در $S_\varphi(M)$ توسیعی دارد به یک $\varphi$-تایپ جهانی که متناهیاً در $M$ محقق می‌شود و تعریف‌پذیر روی $M$ نیز می‌باشد.

ملاحظه ۲.۴.

(۱) اثبات ارائه شده در [۱۶] کاملاً شبیه اثبات اصلی در مقالهٔ گروتندیک (یعنی قضیه ۱.۱) است، البته در زبان نظریهٔ مدل‌ها. همچنین، به گفتهٔ گروتندیک در [۷]، اثبات او مبتنی بر ایده‌ای از اِبرلین[1] است.

(۲) قضیه ۱.۴ قوی‌تر از قضیه ۲.۲ است، زیرا مفهوم «پایداری درون یک مدل» (و نه لزوماً برای یک تئوری) را بررسی می‌کند. در [۲۳]، این دستاورد «قضیهٔ بنیادینِ قویِ نظریهٔ پایداری[2]» نامیده شده است.

(۳) در [۱۶]، پیلی یادآوری می‌کند که در مقالهٔ [۱۵] در دههٔ هشتاد دستاوردی مشابه قضیهٔ گروتندیک را ثابت کرده است (البته در زبان نظریهٔ مدل‌ها) که به‌طور شگفت‌انگیزی اساساً اثباتی مشابه با قضیهٔ گروتندیک دارد. همچنین، یُوینوُ[3] در مقالهٔ [۱۰] قضیهٔ گروتندیک را برای منطق‌های توسیع‌یافته، موسوم به منطق‌های پیوسته، اثبات کرده است.

(۴) به عنوان کاربردی از مفهوم پایداری درون یک مدل، کریوین[4] و مائوری[5] در مقالهٔ [۱۲] دستاورد بسیار مهمی در مورد وجود فضاهای $\ell_p$ ($1 \leq p < \infty$) درون فضاهای باناخِ پایدار ثابت کردند.

(۵) در قضیه فوق، هر تایپ جهانی که توسیع یک تایپ $p \in S_\varphi(M)$ باشد و همچنین در $M$ متناهیاً محقق شود را یک شریکِ ارثِ $p$[6] گویند. این مفهوم از اهمیت بسیاری در نظریه مدل‌ها برخوردار است و کاربردهای مهمی نیز دارد که برای مشاهدهٔ آنها خواننده را به کتاب‌های نظریه مدل‌ها، از جمله [۱۴]، ارجاع می‌دهیم.

این پایان سرگذشت قضیهٔ گروتندیک در نظریهٔ مدل‌ها نیست. در [۱۶]، پیلی مفهومی از یک تایپ «به‌طور ژِنریک پایدار» ارائه می‌دهد که اساساً مناسب برای تئوری‌های «پایدار» است و به‌اندازهٔ‌کافی برای مطالعه در تئوری‌های ناپایدار قوی نیست. او در آنجا، یک $\varphi$-تایپ را به‌طور ژنریک پایدار روی $M$ می‌نامد اگر متناهیاً در $M$ محقق شود و همچنین روی $M$ تعریف‌پذیر باشد. پیلی پرسش می‌کند: «برای چه مفهوم مناسبی روی $M$ (مشابه خاصیت نداشتن ترتیب) می‌توان دستاوردی مشابه قضیهٔ گروتندیک برای مفهوم «قویتری» از به‌طور ژنریک پایدار بودن ارائه داد؟» (دراقع، مفهوم قویترِ به‌طور ژنریک پایدار بودن قبلاً در نظریهٔ مدل‌ها بطور دقیق ارائه شده است.) فصل بعد مطالعه‌ای است برای پاسخ به این پرسش.

---

[1] William Frederick Eberlein (1917-1986)
[2] Strong Fundamental Theorem of Stability (SFTS)
[3] Jose Iovino
[4] Jean-Louis Krivine (1939-)
[5] Bernard Maurey (1948-)
[6] coheir of $p$

## ۵. تایپ‌های به‌طور ژِنریک پایدار، ۲۰۲۱

در نظریهٔ مدل‌ها، مفهومی به نام «تایپ‌های به‌طور ژنریک پایدار» وجود دارد که ریشه‌های آن به «نقطهٔ ژنریک» در کارهای وِی[۱] در هندسه جبری باز می‌گردد. (به زبانی ساده، یک نقطه از یک چندگوناای[۲] $V$ (روی یک میدان $K$) ژنریک نامیده می‌شود اگر در هیچ زیرچندگونای سِرِه (روی $K$) نباشد.) در نظریهٔ مدل‌ها، مفهوم «تایپ‌های به‌طور ژنریک پایدار» را پوآزا[۳]، اساساً برای تئوری‌های پایدار، معرفی کرد. (مراجعه شود به مقالات [۱۸]، [۱۹]، یا کتاب [۲۰].) برای تئوری‌های دلخواه، این مفهوم توسط پیلِی و تَنُویچ[۴] [۱۷] معرفی گردید. در ادامه، این مفهوم معرفی و ارتباط آن با دستاوردهای قبل مطالعه می‌گردد.

تعریف یک «تایپ کامل روی $M$» و «تایپ جهانی» به‌ترتیب مشابه تعریف «$\varphi$-تایپ کامل روی $M$» و «$\varphi$-تایپ جهانی» می‌باشد، با این تفاوت که فقط شامل فرمول $\varphi$ نیستند بلکه شامل همهٔ فرمول‌های زبان $L$ هستند. به‌طور مشابه، مفاهیم «تعریف‌پذیری روی $M$» و «متناهیاً محقق شده در $M$» برای تایپ‌های کامل روی $M$ و تایپ‌های جهانی قابل تعریف‌اند.

فرض کنید $p(x)$ یک تایپ جهانی باشد که در مدل $M$ متناهیاً محقق شود. دنباله‌ای برای $p$ روی $M$ به‌صورت زیر می‌سازیم. عنصر $a_1$ را به‌صورتی انتخاب می‌کنیم که همهٔ فرمول‌های تایپ $p$ که فقط شامل پارامترهایی در $M$ هستند را محقق کند. با استقراء، فرض کنید عناصر $\{a_i : i \leq \alpha\}$ را انتخاب کرده‌ایم. عنصر $a_{\alpha+1}$ را به‌صورتی انتخاب می‌کنیم که همهٔ فرمول‌های تایپ $p$ که فقط شامل پارامترهایی در $M \cup \{a_i : i \leq \alpha\}$ هستند را محقق کند.[۵] دنباله‌ای مانند $\{a_i : i < \alpha\}$ که به این صورت تعریف می‌شود را یک دنبالهٔ مورلی[۶] از $p$ روی $M$ می‌نامیم.[۷] (توجه داریم که دنبالهٔ مورلی لزوماً یکتا نیست، البته هر دو دنباله مورلی تایپ یکسانی دارند. همچنین، هر دنبالهٔ مورلی، تمایزناپذیر[۸] است.)

**تعریف ۱.۵.** فرض کنید $p(x)$ یک تایپ جهانی باشد که در مدل $M$ متناهیاً محقق شود. گوییم $p$ روی $M$ به‌طور ژنریک پایدار است، اگر برای هر دنبالهٔ مورلی $\{a_i : i < \omega + \omega\}$ از $p$ روی $M$، و هر فرمول $\phi(x)$، مجموعه $\{i : \phi(a_i)$ درست است$\}$ متناهی یا متمم‌متناهی باشد.

**ملاحظه ۲.۵.**

---

[1] André Weil (1906-1998)

[2] variety

[3] Bruno Poizat (1946-)

[4] Predrag Tanović (1961-)

[5] توجه داریم که همه این عناصر در مدلِ بزرگ یافت می‌شوند.

[6] Morley sequence

[7] در اینجا، $\alpha$ یک اردینال دلخواه است و می‌تواند بزرگتر از $\omega$ باشد.

[8] indiscernible

(۱) هر تایپ به‌طور ژنریک پایدار (روی $M$) تعریف‌پذیر (روی $M$) می‌باشد. برعکس این مطلب درست نیست، یعنی تایپ‌هایی وجود دارند که همزمان تعریف‌پذیر (روی $M$) و متناهیاً (در $M$) محقق می‌شوند، اما روی هیچ مدلی به‌طور ژنریک پایدار نیستند. (مراجعه شود به [۴].)

(۲) در تئوری‌های $NIP$[۱] هر تایپ که همزمان تعریف‌پذیر (روی $M$) و متناهیاً (در $M$) محقق شود، به‌طور ژنریک پایدار (روی $M$) است. (مراجعه شود به [۹].)

گَنون[۲] [۶]، با استفاده از ایده‌هایی از سیمون[۳] [۲۲]، در قضیهٔ بسیار زیبای ذیل، جنبه‌ای آنالیزی از تایپ‌های به‌طور ژنریک پایدار را روشن کرد.

**قضیه ۳.۵.** فرض کنید $T$ یک تئوری در زبانی شمارا، $M$ مدلی از $T$ و $p(x)$ یک تایپ جهانی باشد. اگر $p$ به‌طور ژنریک پایدار روی $M$ باشد، آنگاه دنباله‌ای مانند $(a_i : i < \omega)$ در $M$ وجود دارد به‌طوری که

$$\lim_i tp(a_i/\mathbb{M}) = p.$$

در قضیهٔ فوق، حد در توپولوژی منطقی (توپولوژی اِستون[۴]) تعریف شده است. به‌طور هم‌ارز، حدِ فوق به این معنا است که، برای هر فرمول $\phi(x)$ با پارامترهایی در مدل جهانی، $\phi(x) \in p$ اگر و تنها اگر سرانجام $\phi(a_i)$ها درست باشند؛ یعنی برای یک عدد طبیعی $n$ داشته باشیم: $\phi(a_i)$ درست است اگر $n \leq i$.

**ملاحظه ۴.۵.** در قضیه ۳.۵، نکتهٔ اساسی وجودِ دنباله‌ای که همگرا به $p$ باشد نیست (زیرا این مطلبی شناخته‌شده است که هر دنبالهٔ مورلی همگرا به $p$ است)، بلکه وجودِ چنین دنباله‌ای در «درون $M$» اهمیت دارد. هرچند که اثباتِ گَنون کاملاً نظریه مدلی است، اما این قضیه ریشه در دستاوردی بسیار مهم در آنالیز تابعی در مقالهٔ بورگِین[۵]، فِرملین[۶]، تالاگراند[۷] [۳] دارد.

پرسشی که به‌طور طبیعی مطرح می‌شود این است که آیا برعکس قضیهٔ فوق برقرار است و آیا ارتباطی بین آن با قضیهٔ گروتندیک و همچنین پرسش پیلِی (در بخش قبل) وجود دارد. در ادامه به این پرسش‌ها پاسخ می‌دهیم.

---

[1] Non Independence Property
[2] Kyle Gannon
[3] Pierre Simon
[4] Stone topology
[5] Jean Bourgain (1954-2018)
[6] David Heaver Fremlin (1942-)
[7] Michel Pierre Talagrand (1952-)

در [۱۱]، برای یک مدل $M$، مفهومی موسوم به «سرانجام وابسته[1]» ارائه شده است که هم‌خانواده مفهوم «وابستگی[2]» در نظریهٔ مدل‌ها می‌باشد. به‌دلیل فنی بودن این مفهوم، خواننده را به تعریف ۱.۳ در آن مقاله ارجاع می‌دهیم. یک مدل $M$ را «سرانجام پایدار» گوییم هرگاه در دو گزارهٔ زیرین صدق کند: (الف) هر فرمول در $M$ پایدار باشد. (مراجعه شود به تعریف ۱.۲) (ب) $M$ سرانجام وابسته باشد.

اکنون آماده‌ایم که دستاوردی مشابه قضیهٔ گروتندیک که مناسب برای تایپ‌های بطور ژنریک پایدار است و همچنین پاسخی به سوال پیلِی می‌باشد را ارائه دهیم. در [۱۱]، با استفاده از ایده‌هایی از [۲۲] و [۶]، قضیهٔ زیر ثابت شده است.

قضیه ۵.۵ (تقویتِ قضیهٔ حدِ دوگانهٔ گروتندیک[3]). فرض کنید $T$ یک تئوری در زبانی شمارا باشد و $M$ مدلی از $T$ باشد. آنگاه، گزاره‌های زیرین هم‌ارز هستند.

(۱) $M$ سرانجام پایدار است.

(۲) هر تایپ کامل روی $M$ توسیعی دارد به یک تایپ جهانی که روی $M$ به‌طور ژنریک پایدار می‌باشد.

ملاحظه ۶.۵.

(۱) قضیهٔ فوق از دو منظر قویتر از قضیه ۱.۴ است. یکم، برای تایپ‌های کامل است و نه فقط برای $\varphi$-تایپ‌ها. توجه داریم که اساساً مفهوم دنبالهٔ مورلی برای $\varphi$-تایپ‌ها قابل‌تعریف نیست. دوم، همان‌طور که قبلاً بیان شد، مفهوم تایپِ به‌طور ژِنریک پایدار قویتر از تایپ‌های تعریف‌پذیر و متناهیاً محقق شده است.

(۲) اثباتِ این قضیه کاملاً با ابزارهای نظریه مدلی امکان‌پذیر است و حتی نیازی به استفاده از قضیه ۱.۱ نیست. (مراجعه شود به قضیه ۱۰.۴ در [۱۱].) همچنین، برای اثبات این قضیه، توصیفِ کاملاً جدیدی از تایپ‌های به‌طور ژِنریک پایدار ارائه شده است ([۱۱، قضیه ۴.۴]) که فی‌نفسه مهم و دارای کاربرد می‌باشد. در واقع، قضیهٔ اخیر کامل کنندهٔ قضیهٔ گنون (قضیه ۳.۵) است. برای اثبات آن، مفهومی از همگرایی استفاده می‌شود که قوی‌تر از همگرایی در قضیهٔ گنون است، که به‌دلیل فنی بودن این مفهوم خواننده را به [۱۱، تعریف ۴.۳] ارجاع می‌دهیم.

(۳) قضیه ۵.۵ را می‌توان به زبانِ آنالیز تابعی ترجمه کرد و دستاوری مشابه قضیه ۱.۱ به‌دست آورد. محتمل است که چنین دستاوردی کاربردهایی در آینده داشته باشد.

[1] eventually dependence
[2] Dependence property
[3] Reinforcement of Grothendieck's double limit theorem (RDLT)

## ۶. سخن پایانی

گرچه سرگذشت قضیهٔ گروتندیک در نظریهٔ مدل‌ها مسیر نسبتاً طولانی طی نموده است، ولی این سرگذشت می‌تواند همچنان ادامه داشته باشد. پرسش‌های بسیاری پیرامون آن می‌توانند مطرح گردند، که یکی از آنها تعمیم قضیه ۵.۵ برای اندازه‌های کِزْلِر[1] است. اندازه‌های کزلر تعمیم طبیعی مفهوم تایپ‌ها هستند که کاربردهای بسیار مهمی در نظریهٔ مدل‌ها داشته‌اند (مراجعه شود به [۸] و [۹]). به‌طور طبیعی، ارائه یک مفهوم مناسب برای اندازه‌های کزلرِ به‌طور ژنریک پایدار از اهمیت ویژه‌ای برخوردار است. همچنین، ارائه مفهوم مناسب و اثبات یک هم‌ارزی مشابه با قضیه ۵.۵ می‌تواند سرگذشتی دیگر در مقاله‌ای دیگر باشد.

اجازه دهید که در انتها دیدگاهِ فلسفهِ ریاضیِ این مقاله را تصریح کنیم. در مقالهٔ بسیار زیبای [۱]، مایکل عطیه[2] تأکید می‌کند که یکی از مهمترین اهداف در ریاضیات می‌بایست «متّحد کردنِ ریاضیات» باشد. او بحث می‌کند که چرا عمیق‌ترین، زیباترین، و کاربردی‌ترین بخش‌های ریاضیات آنهایی هستند که ارتباطی بین حوزه‌های مختلف ایجاد می‌کنند. در واقع، عطیه عمر خود را صرفِ متّحد کردنِ شاخه‌هایِ مختلف ریاضیات، و حتی در کارهای متأخر خود، متّحد کردن «فیزیک» و «ریاضیات» کرد. قطعاً، سرگذشت قضیهٔ حد دوگانهٔ گروتندیک یکی از مصادیق تأکید عطیه بر «متّحد کردنِ ریاضیات» است.

## تشکر و قدردانی

در اینجا قصد دارم که از جناب آقای دکتر محسن خانی به دلیل مطالعه نسخه اولیه این مقاله و پیشنهادات بسیار مفید ایشان تشکر و قدردانی نمایم. از پژوهشگاه دانش‌های بنیادین (IPM) تشکر می‌کنم. این اثر تحت حمایت مالی IPM با شماره طرح «۱۴۰۱۰۳۰۱۱۷» انجام شده است.

## مراجع

[1] Atiyah, M., (1978) The unity of mathematics, *Bulletin of the London Mathematical Society*, vol. 10, pp. 69–76.

[2] Ben Yaacov, I., (2014) Model theoretic stability and definability of types, after A. Grothendiek, *Bulletin of Symbolic Logic*, 20, pp 491-496.

[3] Bourgain, J., Fremlin, D.H., and Talagrand, M., (1978) Pointwise compact sets of Baire-measurable functions. *American Journal of Mathematics*, 100(4):pp. 845-886 1978.

[1] Keisler measures
[2] Michael Atiyah (1929-2019)

[4] Conant, G., Gannon, K., Hanson, J., (2023) Keisler measures in the wild, *Journal of Model Theory*, to appear.

[5] Fremlin, D., *Measure Theory, vol.4, Topological Measure Spaces*, Torres Fremlin, Colchester, (2006).

[6] Gannon, K., (2022) Sequential approximations for types and Keisler measures, *Fundamenta Mathematicae*, 257, 305-336.

[7] Grothendieck, A., (1952) Critères de compacité dans les espaces fonctionnels généraux, *American Journal of Mathematics* 74, 168-186.

[8] Hrushovski, E., Peterzil, Y., and Pillay, A., (2008) Groups, measures, and the NIP, *Journal of the American Mathematical Society*, vol. 21, no. 2, pp. 563–596.

[9] Hrushovski, E., and Pillay, A., (2011) On NIP and invariant measures, *Journal of the European Mathematical Society*, 13, 1005-1061.

[10] Iovino, J., (1999) Stable models and reflexive Banach spaces. *Journal of Symbolic Logic*, 64, 1595-1600.

[11] Khanaki, K., (2022) Remarks on convergence of Morley sequences, https://arxiv.org/abs/2110.15411

[12] Krivine, J.-L., and Maurey, B., (1981) Espaces de Banach stables. *Israel J. Math.*, 39(4):273-295.

[13] Pietsch, P., *History of Banach Spaces and Linear Operators*. Birkhäuser Boston, MA, (2008).

[14] Pillay, A., *An Introduction to Stability Theory*, Oxford Logic Guides (1983).

[15] Pillay, A., (1982) Dimension and homogeneity for elementary extensions of a model. *The Journal of Symbolic Logic*, 47(1): 147-160.

[16] Pillay, A., (2016) Generic Stability and Grothendieck, *South American Journal of Logic*, Vol. 2, n. 2, pp. 437-442

[17] Pillay, A., and Tanović, P., (2011) Generic stability, regularity, and quasiminimality, Models, logics, and higher-dimensional categories, *CRM Proc. Lecture Notes, vol. 53, Amer. Math. Soc., Providence, RI*, pp. 189–211. MR 2867971

[18] Poizat, B., (1983) Groupes stables avec types génériques réguliers, *Journal of Symbolic Logic*, Vol. 48, No. 2 (Jun., 1983), pp. 339-355

[19] Poizat, B., (1981) Sous-groupes définissables d'un groupe stable, *Journal of Symbolic Logic*, Vol. 46, No. 1 (Mar., 1981), pp. 137-146.

[20] Poizat, B., *Stable groups*, Mathematical Surveys and Monographs, vol. 87, Providence, RI: American Mathematical Society, pp. xiv+129, (2001) doi:10.1090/surv/087, ISBN 0-8218-2685-9, MR 1827833 (Translated from the 1987 French original.)

[21] Shelah, S., (1971) Stability, the f.c.p., and superstability; model theoretic properties of formulas in first order theory, *Annals of Mathematical Logic*, vol. 3, no. 3, pp. 271-362.

[22] Simon, P., (2015) Invariant types in NIP theories, *Journal of Mathematical Logic*, Vol. 15, No. 02, 1550006.

[23] Starchenko, S., (2017) On Grothendieck's approach to stability, unpublished note, https://www3.nd.edu/ sstarche/papers/groth-stab.pdf

[24] Young, N.J., (1971) On Pták's double-limit theorems, *Proceedings of the Edinburgh Mathematical Society*, Volume 17 , Issue 3 , pp. 193 - 200

(کریم خانکی) دانشگاه صنعتی اراک، دانشکده علوم پایه و پژوهشگاه دانش‌های بنیادین (IPM)

https://www.arakut.ac.ir/fa/khanaki.html :

khanaki@arakut.ac.ir :